\input amstex
\loadeufm
\documentstyle{amsppt}
\magnification=1200
%\hcorrection{.25in}
%\advance\vsize-.75in
\TagsOnRight

\define\countable{\operatorname{countable}}

\define\ddef{\operatorname{def}}

\define\End{\operatorname{End}}

\define\GG{\operatorname{G}}

\define\GL{\operatorname{GL}}

\define\Hg{\operatorname{Hg}}

\define\Hom{\operatorname{Hom}}

\define\Lie{\operatorname{Lie}}

\define\Spec{\operatorname{Spec}}

% Bourbaki letters
\define\C{\Bbb C} % complex numbers
 % quaternions
 % natural numbers
\define\Q{\Bbb Q} % rational numbers
 % real numbers
\define\Z{\Bbb Z} % integers

\topmatter
\title
On the Hodge group and invariant cycles
\endtitle

\title
{\bf On the Hodge group and invariant cycles of a simple Abelian variety with a stable reduction of odd toric rank}
\endtitle
\footnotetext"\ "{This work was carried out with the financial support of
the Russian Foundation for
Basic Research
(grant N 18-01-00143).\newline
}

\author
O.V. Oreshkina 
\endauthor

\author
O.V. Oreshkina 
\endauthor

\address  Vladimir State University
\endaddress

\email
papichonok\@yandex.ru
\endemail

\abstract
The semi-simplicity of the Hodge group is proved for a simple Abelian variety with a stable reduction of odd toric (reductive) rank.
If, besides, the dimension of the Abelian variety is an odd integer, then the Hodge conjecture on algebraic cycles holds for it.

The bibliography: 13 titles.


\endabstract

\endtopmatter

\document

\head
{\bf Introduction }
\endhead

It is well known that the Hodge group of a complex Abelian variety with a totally degenerated reduction is semi-simple (\cite{1}, Theorem 4.1).

By definition, the toric rank of the special fibre ${\Cal M}_v$ of the minimal N\'eron model of an Abelian variety over the fraction field of a discrete valuation ring ${\Cal O}_v$ with a stable reduction at the closed place $v \in \Spec {\Cal O}_v$ is equal to the dimension of the greatest linear subtorus of the reductive algebraic group $\Cal M_v$ (this rank is called also a reductive rank (\cite{2}, Section 2.1.12).

If the toric rank $r_v$ of an absolutely simple Abelian variety $J$ is an odd integer, then this elementary condition allows us to prove the semi-simplicity of the Hodge group of the Abelian variety $J$ and, if the dimension $g = \dim J$ is an odd integer, then the Hodge conjecture on algebraic cycles holds for the variety $J$ (Theorem 1).
If, besides, the Abelian variety $J$ is the generic scheme fibre of an Abelian scheme $\pi_{X / C} : X \to C$ over some affine curve $C$ with a stable reduction of odd toric rank at some infinite place of the curve $C$ then, for any Abelian scheme $\pi_{Y / C} : Y \to C$, there is a canonical isomorphism
$\Hom_C(X, Y) \,\widetilde{\rightarrow}\, \Hom(R_1 \pi_{X / C\ast} \Z, R_1 \pi_{Y / C\ast} \Z)$ and,
for the Abelian scheme $\pi_{X / C} : X \rightarrow C$, the Grothendieck conjecture on invariant cycles holds (Theorem 2).

\head
{\bf Formulations and proofs of basic results}
\endhead

{\bf Theorem 1.} {\it Let $J$ be an absolutely simple $g$-dimensional Abelian variety over a subfield $k \subset \C, \,\End_{\C}(J) = \End_k(J), \, v$ a discrete valuation of the field $k$. If the Abelian variety $J$ has a stable reduction at the place $v$ and the toric rank $r_v$ of the special fibre $\Cal{M}_v$ of the minimal N\'eron model is an odd integer, then the Hodge group $\Hg(J)$ of the Abelian variety $J$ is semi-simple, the division $\Q$-algebra $\End_{\C}(J) \otimes _{\Z} \Q$ is a totally real field of odd degree $e$ over the field $\Q$. If, besides, the number $g$ is odd, then the semi-simple Lie algebra $\Lie \Hg(J)\otimes_{\Q} \C$ has type
$C_{g/e} \times \cdots \times C_{g/e} = C^{\times e}_{g/e}$ and, for any integer $p \in \{1, \ldots , g\}$, the $\Q$-space of Hodge cycles
$$[\wedge^{2p}H^1(J, \Q)]^{\Hg(J)}=H^{2p}(J, \Q) \cap H^{p,p}(J, \C)$$
is generated by classes of intersections of divisors on the Abelian variety $J$ $($in particular, for the variety $J$, the Hodge conjecture on algebraic cycles holds$)$.}

{\it Proof}. We denote the fraction field of the ring $\Cal{O}_v$ of the
discrete valuation by $k_v$. Since a reduction is stable, then the fibre of the minimal N\'eron model $\Cal{M} \rightarrow \Spec \Cal{O}_v$ over the closed point $v \in \Spec \Cal O_v$ is an extension of an Abelian variety by a linear torus of dimension $r_v$. It is well known that the formation of N\'eron models commutes with \'etale base change. Therefore, considering a suitable base change, which is defined by the finite \'etale morphism
$\Spec \widetilde{\Cal O}_v \rightarrow \Spec \Cal O_v$, we may assume that a linear torus is totally decomposed over the residue field $\kappa(v)$ of the place $v$.

Denoting by $\Cal M^0_v$ the connected component of the neutral element of the algebraic group $\Cal M_v$, we have an exact sequence of algebraic groups over the field $\kappa(v)$
$$
1 \rightarrow \GG^{r_v}_m \rightarrow \Cal M^0_v \rightarrow A \rightarrow 0,
\tag{1}
$$
where $A$ is some Abelian variety over the field  $\kappa(v)$. It is well known that the exact sequence \thetag{1} yields an extension

$$1 \rightarrow [\GG_m \otimes_{\kappa(v)} \overline{\kappa(v)}]^{r_v}
\rightarrow \Cal M^0_v \otimes_{\kappa(v)} \overline{\kappa(v)}
\rightarrow A \otimes_{\kappa(v)} \overline{\kappa(v)} \rightarrow 0$$
of algebraic groups over the field $\overline{\kappa(v)}$ which is given by a point of the variety
$[A^{\vee} \otimes_{\kappa(v)} \overline{\kappa(v)}]^{r_v}$, where $A^{\vee}$ is an Abelian variety dual to the Abelian variety $A$ (\cite{3}, Ch. VII, \S \, 3, Section 16, Comment to Theorem 6).

In virtue of the universal property of the N\'eron model (\cite{2}, Section (1.1.2)), there is the canonical isomorphism
$$\End_{\Spec \Cal O_v}(\Cal M) \,\,\widetilde{\rightarrow}\,\, \End_{k_v}(J).$$

Therefore, in virtue of the well-known equalities
$$
\Hom_{\kappa(v)}(\GG_m, A) = \Hom_{\kappa(v)}(A, \GG_m) = 0
\tag{2}
$$
in the category of commutative algebraic groups over the field $\kappa(v)$, the canonical morphisms of rings
$$
\End_k(J) \,\widetilde{\rightarrow}\, \End_{k_v}(J) \,\widetilde{\rightarrow}\,
\End_{\Spec \Cal O_v}(\Cal M) \rightarrow \End_{\kappa(v)}(\Cal M^0_v) \rightarrow \End_{\kappa(v)}(\GG^{r_v}_m)
\tag{3}
$$
are defined
because any endomorphism $\varphi \in \End_{\Spec \Cal O_v}(\Cal M)$ defines an endomorphism of the special fibre $\varphi_v \in \End_{\kappa(v)}(\Cal M_v)$, which in turn gives an endomorphism of the connected component of the neutral element of the special fibre $\varphi^0_v \in \End_{\kappa(v)}(\Cal M^0_v)$; clearly the algebraic group $\varphi^0_v (\GG^{r_v}_m)$ is a linear subtorus in the group $\Cal M^0_v$, so that the formula \thetag{2} gives the inclusion $\varphi^0_v (\GG^{r_v}_m) \subset \GG^{r_v}_m$.

It is obvious that the image of the endomorphism
$J @>{n}>> J$ of the multiplication by the integer $n \geq 2$
in the ring $\End_{\kappa(v)}(\GG^{r_v}_m) = M_{r_v}(\Z)$ is non-trivial and it is represented in the ring of matrices
$M_{r_v}(\Z)$ by a scalar matrix of homothety with the coefficient $n$. Therefore, from the simplicity of the  Abelian variety $J$, it follows that there is a canonical embedding of the division $\Q$-algebra
$E \,\,\overset{\ddef}\to=\,\, \End_k(J) \otimes_{\Z} \Q$ into the ring $M_{r_v}(\Q)$, which defines the structure of a left $E$-module on the $\Q$-space $\Q^{\oplus r_v}$.

Assume that $\dim_{\Q} E$ is an even integer. Then, in virtue of the oddness of the integer $r_v$, in the $\Q$-space $\Q^{\oplus r_v}$ there exists a 1-dimensional subspace $L$, which is annihilated by the division $\Q$-algebra $E$. In this case, the integer $n \in E$ annihilates the line $L$, but this is impossible, because on the space $L$ this integer induces the homothety with coefficient $n$. Therefore, the $\Q$-algebra $\End_k(J) \otimes_{\Z} \Q$ has an odd dimension over the field $\Q$. The Albert classification of division $\Q$-algebras $\End_{\C}(J) \otimes_{\Z} \Q$
(\cite{4}, Ch. IV, \S \, 21, Theorem 2) shows that the $\Q$-algebra $\End_{\C}(J) \otimes_{\Z} \Q$
is a totally real field of {\it odd} degree $e$ over the field $\Q$, so that the Hodge group is semi-simple
(\cite{5}, Lemma 1.4). It is well known that the canonical representation of the Lie algebra $\Lie \Hg(J) \otimes_{\Q} \C$ in the space $H^1 (J, \C)$ is determined by minuscule weights
(\cite{6}; \cite{7}, Theorem 0.5.1), therefore, if the integer $g$ is odd, then the semi-simple Lie algebra $\Lie \Hg(J) \otimes_{\Q} \C$ has type $C_{g/e} \times \cdots \times C_{g/e} = C^{\times e}_{g/e}$ and, for any integer $p \in \{1, \ldots , g\}$, the $\Q$-space of {\it Hodge cycles}
$$[\wedge^{2p}H^1(J, \Q)]^{\Hg(J)}=H^{2p}(J, \Q) \cap H^{p,p}(J, \C)$$
is generated by classes of intersections of divisors on the Abelian variety $J$ (\cite{8}, Theorem 5.1). Theorem is proved.

{\bf Corollary 1.} {\it Let $J$ be an absolutely simple $g$-dimensional Abelian variety over a subfield $k \subset \C$, $\End_k(J) = \End_{\C}(J)$, \,\, \, $v$ a discrete valuation of the field $k$.
If the Abelian variety $J$ has a stable reduction at the place $v$ and the toric rank of the special fibre $\Cal M_v$ of the minimal N\'eron model is equal to $1$, then $\End_{\C}(J) = \Z$ and the Hodge group $\Hg(J)$ is $\Q$-simple.}

{\it Proof}. In the case under consideration, $r_v = 1$, therefore from \thetag{3} it follows that there exists a canonical non-trivial morphism of division $\Q$-algebras
$$\End_k(J) \otimes_{\Z} \Q \rightarrow \End_{\kappa(v)}(\GG_m) \otimes_{\Z} \Q = \Q.$$
Therefore $\End_k(J) = \Z$. Consequently, the Hodge group $\Hg(J)$ is $\Q$-simple by M.V. Borovoi's theorem  \cite{9}.

The Grothendieck conjecture on invariant cycles (\cite{10}, P. 103; \cite{11}, Conjecture B.181; \cite{12}, P. 26) states that, if $\pi_{X / C} : X \rightarrow C$ a is smooth projective morphism of smooth connected quasi-projective complex varieties and $\alpha \in H^0 (C, R^{2p}{\pi_{X / C^{\ast}}} \Q)$
(where $p$ is a natural integer) so that, for some point $s_0 \in C$, the restriction
$\alpha \vert_{X_{s_0}} \in H^{2p} (X_{s_0}, \Q)$ is an algebraic cohomology class then, for any point $s \in C$, the class $\alpha \vert_{X_s} \in H^{2p} (X_s, \Q)$ is algebraic.

{\bf Theorem 2.} {\it Let $C$ be a smooth complex affine curve, $\overline{C}$ its smooth projective model, $\pi_{X / C} : X \rightarrow C$ an Abelian scheme of odd relative dimension so that the generic scheme fibre of the structure morphism $\pi_{X / C} $ is an absolutely simple Abelian variety.
Assume that the fibre $\Cal M_v$ of the minimal N\'eron model $\Cal M \rightarrow \overline{C}$ of the Abelian scheme $\pi_{X / C} : X \rightarrow C$ over some infinite place $v \in \overline{C} \setminus C$ is an extension of an Abelian variety by a linear torus of odd dimension.

Then, for any Abelian scheme $\pi_{Y / C} : Y \rightarrow C$, there is the canonical isomorphism
$$\Hom_C(X, Y) \,\widetilde{\rightarrow}\, \Hom(R_1 \pi_{X / C{\ast}} \Z, R_1 \pi_{Y / C{\ast}} \Z)$$
and, for any natural integer $p$, the space of invariant cycles
$$H^0 (C, R^{2p}\pi_{X / C{\ast}} \Q) \,\widetilde{\rightarrow}\, H^{2p}(X_s, \Q)^{\pi_1(C,s)}$$
is generated by cohomology classes of algebraic cycles on the Abelian variety $X_s$. In particular,  the Grothendieck conjecture on invariant cycles holds for the Abelian scheme $\pi_{X / C} : X \rightarrow C$.}

{\it Proof}. Let $\eta$ be the generic point of the curve $\overline{C}$, \,\,\, $\widetilde{\overline{C}} \rightarrow \overline{C}$ a ramified covering,  $\widetilde{\eta}$ the generic point of the curve
$\widetilde{\overline{C}}$, \,\,\, $\widetilde{\Cal M} \rightarrow \widetilde{\overline{C}}$ the minimal N\'eron model of the Abelian variety
$X_{\widetilde{\eta}}= X_{\eta} \otimes_{\kappa(\eta)} \kappa(\widetilde{\eta})$.
If a reduction is stable, then the connected component of the neutral element of the special fibre of the N\'eron model $\Cal M \rightarrow {\overline{C}}$ is isomorphic to the connected component of the neutral element of an appropriate special fibre of the N\'eron model $\widetilde{\Cal M} \rightarrow \widetilde{\overline{C}}$ (\cite{2}, Corollary 3.3, Corollary 3.9). In particular, a toric rank of corresponding special fibres
is preserved under the base change, which is defined by a ramified coverings
$\widetilde{\overline{C}} \rightarrow \overline{C}$.

On the other hand, the toric rank at the place $v$ is non-trivial, therefore, in virtue of the absolute simplicity of the Abelian variety $X_{\eta}$, the Abelian scheme $\pi_{X / C} : X \rightarrow C$ has a trivial trace (in other words, in the case under consideration, it is non-isotrivial, for any surjective morphism $S \rightarrow C$ the Abelian scheme $X \times_C S \rightarrow S$ is non-constant).

Using (if necessary) the base change, which is defined by a finite ramified covering $S \rightarrow C$, we obtain from Theorem 1 that the division $\Q$-algebra
$$\End_{\kappa(\overline{\eta})}(X_{\overline{\eta}}) \otimes_{\Z} \Q$$ is a totally real field of odd degree over the field $\Q$.
Therefore, in virtue of the G.A. Mustafin's theorem (\cite{13}, Theorem 4.1), for any Abelian scheme
$\pi_{Y / C} : Y \rightarrow C$, there is the canonical isomorphism
$$\Hom_C(X, Y) \,\widetilde{\rightarrow}\, \Hom(R_1 \pi_{X / C{\ast}} \Z, R_1 \pi_{Y / C{\ast}} \Z).$$

By Deligne's theorem (\cite{7}, Theorem 7.3) there exists a countable subset of closed points
$\Delta_{\countable} \subset C$, such that, for any point $s \in C \setminus \Delta_{\countable}$, the connected component of the unity $G^0$ of the closure $G$ of the image of the monodromy representation
$\pi_1(C, s) \rightarrow \GL(H^1 (X_s, \Q))$ in the $\Q$-Zariski topology of the group $\GL(H^1 (X_s, \Q))$ is a normal subgroup of the Hodge group $\Hg(X_s)$ of the Abelian variety $X_s$. In the case under consideration, the Hodge group $\Hg(X_{\overline{\eta}})$ of the generic geometric fibre $X_{\overline{\eta}}$ is a $\Q$-simple algebraic group (\cite{13}, Derivation of Theorem 4.1 from Lemmas 1 - 3). Therefore, from the non-isotriviality of the Abelian scheme $\pi_{X / C} : X \rightarrow C$ it follows that $G^0 = \Hg(X_s)$ for any point $s \in C \setminus \Delta_{\countable}$ (\cite{13}, Proposition 4.1).
According to Theorem 1, the $\Q$-space of Hodge cycles
$$[\wedge^{2p}H^1(X_{\overline{\eta}}, \Q)]^{\Hg(X_{\overline{\eta}})}
= H^{2p}(X_{\overline{\eta}}, \Q) \cap H^{p,p}(X_{\overline{\eta}}, \C)$$
is generated by classes of intersections of divisors on the Abelian variety $X_{\overline{\eta}}$. Therefore,
from the existence of obvious isomorphisms and the embedding
$$H^0 (C, R^{2p}\pi_{X / C\ast} \Q) \,\widetilde{\rightarrow}\, H^{2p} (X_s, \Q)^{\pi_1(C, s)}
\hookrightarrow H^{2p}(X_s, \Q)^{G^0} \,\widetilde{\rightarrow}\, H^{2p}(X_{\overline{\eta}}, \Q)^{G^0}$$
$$= [\wedge^{2p}H^1(X_{\overline{\eta}}, \Q)]^{\Hg(X_{\overline{\eta}})}
= H^{2p}(X_{\overline{\eta}}, \Q) \cap H^{p,p}(X_{\overline{\eta}}, \C),$$
it follows that the Grothendieck conjecture on invariant cycles holds for the Abelian scheme $\pi_{X / C} : X \rightarrow C$. Theorem is proved.

The author is grateful to S.G. Tankeev for an attention to the article.

\enddocument